\documentclass[secthm,seceqn,amsthm,ussrhead,10pt]{amsart}
\usepackage{amsmath,latexsym}
\usepackage[english]{babel}
\usepackage[psamsfonts]{amssymb}
\usepackage{times}
\usepackage{cite}

\usepackage[mathcal]{euscript}

\newtheorem{Th}{Theorem}
\newtheorem{Lem}[Th]{Lemma}

\newtheorem{df}[Th]{Definition}

\begin{document}
	\sloppy
	
{\huge Generalized derivations of multiplicative $n$-ary $Hom$-$\Omega$ color algebras}
\footnote{
The authors were supported by Funda\c{c}\~{a}o para a Ci\^{e}ncia e a Tecnologia (Portugal), project PEst-OE/MAT/UI0212/2015 of CMA-UBI and
by Ministerio de Econom\'{i}ıa y Competitividad (Spain), project MTM2013-45588-C3-1-P and by RFBR 16-31-00096. }
\medskip

\medskip
\textbf{P. D. Beites$^{a}$, Ivan Kaygorodov$^{b}$, Yury Popov$^{c,d}$}

\medskip

{\tiny
$^a$ Universidade da Beira Interior, Covilh\~{a}, Portugal

$^b$ Universidade Federal do ABC, CMCC, Santo Andr\'{e}, Brasil

$^c$ Universidade Estadual de Campinas, IMECC-UNICAMP, Campinas, Brasil

$^d$ Novosibirsk State University, Novosibirsk, Russia

\medskip

    E-mail addresses:

    Patr\'{i}cia Damas Beites (pbeites@ubi.pt),

    Ivan Kaygorodov (kaygorodov.ivan@gmail.com),
    
    Yury Popov (yuri.ppv@gmail.com).
    
\medskip
    
}

{\bf Abstract.}
We generalize the results of Leger and Luks, 
 Zhang R. and Zhang Y.; Chen, Ma, Ni, Niu, Zhou and Fan; Kaygorodov and Popov (see, \cite{LL00,ZZ10,cml13,cml14, cml15, ZCM16, Z15, ZF16, ZF16col,KP16}) about generalized derivations of 
color $n$-ary algebras to the case of  $n$-ary $Hom$-$\Omega$ color algebras.
Particularly, we prove some properties of generalized derivations of multiplicative  $n$-ary $Hom$-$\Omega$ color algebras.
 Moreover, we prove that the quasiderivation algebra  of any  multiplicative  $n$-ary $Hom$-$\Omega$ color algebra can be embedded into
the derivation algebra of a larger multiplicative  $n$-ary $Hom$-$\Omega$ color algebra.

\


\section*{Introduction}

It is well known that the algebras of derivations and generalized derivations are very important in the study of Lie algebras and its generalizations.
There are many generalizations of derivations
(for example, Leibniz derivations \cite{KP} and Jordan derivations \cite{her}).
The notion of a $\delta$-derivation appeared in the paper of Filippov \cite{fil1},
in which he studied $\delta$-derivations of prime Lie and Malcev algebras \cite{fil2,fil3}.
After that, $\delta$-derivations of
Jordan and Lie superalgebras were studied in \cite{kay1, kay_lie, kay_lie2,kay12mz, zhel, Kokh, zus10} and many other works.
The notion of a  generalized derivation is a generalization of $\delta$-derivation.
The most important and systematic research on the generalized derivations algebras of a Lie algebra and their subalgebras was
due to Leger and Luks \cite{LL00}.
In their article, they studied properties of generalized derivation algebras and their subalgebras, for example, the quasiderivation algebras. They have determined the structure of algebras of quasiderivations and generalized derivations and proved that the     quasiderivation algebra of a Lie algebra can be embedded into the derivation algebra of a larger Lie algebra.
Their results were generalized by many authors. For example, 
Zhang and Zhang \cite{ZZ10} generalized the above results to the case of Lie superalgebras; Chen, Ma, Ni, Zhou and Fan considered  generalized derivations of 
color Lie algebras, 
$Hom$-Lie superalgebras, 
Lie triple systems, $Hom$-Lie triple systems 
and n-$Hom$ Lie superalgebras \cite{cml13,cml14,cml15, Z15, ZCM16,ZF16,ZF16col,KP16com}.
Generalized derivations of color $n$-ary $\Omega$-algebras were studied in \cite{KP16}.

Generalized derivations of simple algebras and superalgebras
were investigated in \cite{KN03,lesha12,lesha14,GP03}.
P\'{e}rez-Izquierdo and Jim\'{e}nez-Gestal used  generalized derivations to study non-associative algebras \cite{GP08,GP09}.
Derivations and generalized derivations of $n$-ary algebras were considered in \cite{bn16,poj06,ck10,kay12izv,kay11aa,kay14sp,kay14mz,Williams} and other works.
For example, Williams proved that, unlike the case of binary algebras, for any $n \geq 3$ there exist a non-nilpotent $n$-Lie algebra with invertible derivation  \cite{Williams},
Kaygorodov described $(n+1)$-ary derivations of simple $n$-ary Malcev algebras \cite{kay14sp} and  generalized derivations algebras of semisimple Filippov algebras over an algebraically closed field of characteristic zero \cite{kay14mz}.

The study of $Hom$-structures was started in the classic paper of Hartwig, Larsson and Silvestrov \cite{hom06}.
After $Hom$-Lie algebras were studied some questions related to 
$Hom$-Lie bialgebras \cite{hombi},
$Hom$-Lie color algebras \cite{homliecol},
$Hom$-alternative and $Hom$-Malcev algebras \cite{homalt},
ternary $Hom$-Nambu-Lie algebras \cite{ternhom},
$Hom$-Sabinin algebras \cite{cm15} and many others.

The main purpose of our work is to generalize the results of 
Leger and Luks \cite{LL00}; Zhang R. and Zhang Y. \cite{ZZ10}; Chen, Ma, Ni, Niu, Zhou and Fan \cite{cml13,cml14, cml15, ZCM16, Z15, ZF16, ZF16col}; Kaygorodov and Popov  \cite{KP16}  to the case of multiplicative $n$-ary $Hom$-$\Omega$ color algebras for an arbitrary variety of $Hom$-identities $Hom$-$\Omega.$
Particularly, we prove some properties of generalized derivations of multiplicative $n$-ary $Hom$-$\Omega$ color algebras;
prove  that the quasiderivation algebra of a multiplicative $n$-ary $Hom$-$\Omega$ color algebra can be embedded into
the derivation algebra of a larger multiplicative $n$-ary $Hom$-$\Omega$ color algebra.

\section{Preliminaries}

In this section we consider some well known $Hom$-algebraic structures that are twisted versions of the original algebraic structures. In fact, they are ($n$-ary) algebras where the identities defining the structure are twisted by homomorphisms, usually called twisting maps.
All of the following $Hom$-algebras are particular cases of
$n$-ary $Hom$-$\Omega$ color algebras that we consider in the present work.

\begin{df}{
Let $\mathbb{F}$ be a field and $\mathbb{G}$ be an abelian group.
A map $\epsilon: \mathbb{G} \times \mathbb{G} \rightarrow \mathbb{F}^*$ is called a bicharacter on $\mathbb{G}$ if the
following relations hold for all $f,g,h \in \mathbb{G}:$

(1) $\epsilon(f,g+h)=\epsilon(f,g)\epsilon (f,h);$

(2) $\epsilon(g+h,f)=\epsilon(g,f)\epsilon(h,f);$

(3) $\epsilon(g,h)\epsilon(h,g)=1.$}
\end{df}

\begin{df}{
A color  $n$-ary algebra $T$ is an $n$-ary $\mathbb{G}$-graded vector space $T=\bigoplus_{g\in \mathbb{G}} T_g$ with a graded $n$-linear map
$[\cdot , \ldots, \cdot ]: T \times \ldots \times T \rightarrow T$  satisfying
$$[T_{\theta_1},  \ldots, T_{\theta_n}] \subseteq T_{\theta_1+ \ldots +\theta_n}, \ \theta_i \in \mathbb{G}.$$}
\end{df}

The main examples of color $n$-ary algebras
are color Lie algebras \cite{cml13},
color Leibniz algebras \cite{dzhuma},
Filippov ($n$-Lie) superalgebras \cite{poj08patricia,ck10,poj03,poj09,poj08}
and $3$-Lie color algebras \cite{ZT15}.

Let $T=\bigoplus_{g\in \mathbb{G}} T_g$ be a color algebra.
An element $x$ is called a homogeneous element of degree $t \in \mathbb{G}$ if $x \in T_t$. We denote this by $hg(x)=t$.
A linear map $D$ is homogeneous of degree $t$ if $D(T_g) \subseteq T_{g+t}$ for all $g \in \mathbb{G}.$
We denote this by $hg(D)=t$. If $t=0$ then $D$ is said to be even. From now on, unless stated otherwise, we assume that all elements and maps are homogeneous. Let $\epsilon$ be a bicharacter on $\mathbb{G}$. For two homogeneous elements $a$ and $b$ we set $\epsilon(a,b):=\epsilon(hg(a),hg(b)).$

Let $D_1, D_2$ be linear maps on $T$ and $x_1, \ldots, x_n \in T.$ We use the following notations:
\begin{eqnarray}\label{colorcom}[D_1,D_2] = D_1D_2 - \epsilon(D_1,D_2)D_2D_1,
\end{eqnarray}
\begin{eqnarray*}X_i = hg(x_1)+\ldots +hg(x_i).
\end{eqnarray*} In particular, we set $X_0=0$.

By ${\rm End}(T)$ we denote the set of all linear maps of
$T$. Obviously, ${\rm End}(T) = \bigoplus_{g \in \mathbb{G}}{\rm End}_g(T)$ endowed with the color bracket (\ref{colorcom}) is a color Lie algebra over $\mathbb{F}$.

\subsection{Binary $Hom$-algebras.}
In this subsection we give some definitions of $Hom$-algebras with a binary multiplication.
It is known that for any family of polynomial identities $\Omega$,
every $\Omega$-algebra is a $\Omega$-superalgebra
and every $\Omega$-superalgebra is a  color $\Omega$-algebra.
 In what follows we only give the definitions of color $Hom$-$\Omega$ algebras.
The definitions of $Hom$-$\Omega$ superalgebras and $Hom$-$\Omega$ algebras 
 are particular cases of the former ones for $\mathbb{G}=\mathbb{Z}_2$ and $\mathbb{G}=\{ 0 \},$ respectively.
On the other hand, the definitions of 
$\Omega$-algebras, $\Omega$-superalgebras and color $\Omega$-algebrascan be obtained from suitable  $Hom$-definitions for $\alpha= id.$

\begin{df}{A color $Hom$-Lie algebra $(T,[\cdot,\cdot], \epsilon, \alpha)$ is a $\mathbb{G}$-graded vector space $T=\bigoplus_{g \in \mathbb{G}} T_g$ with a bicharacter $\epsilon$, an even bilinear map $[\cdot,\cdot]: T \times T \rightarrow T$ and an even linear map $\alpha: T \rightarrow T$ satisfying
$$[x,y]=-\epsilon(x,y)[y,x],$$
$$\epsilon(z,x)[\alpha(x),[y,z]] + \epsilon(x,y)[\alpha(y),[z,x]]
+ \epsilon(y,z)[\alpha(z),[x,y]]=0.$$
}\end{df}

\begin{df}{
A color $Hom$-Jordan algebra $(T,[\cdot,\cdot], \epsilon, \alpha)$ is a 
$\mathbb{G}$-graded vector space  $T=\bigoplus_{g \in \mathbb{G}} T_g$
with a bicharacter $\epsilon$, an even  bilinear map $ [\cdot,\cdot]: T \times T \rightarrow T$ and an even linear map $\alpha: T \rightarrow T$ satisfying
$$[x,y] = \epsilon(x,y) [y,x],$$
$$\epsilon(z, x+w)as_{T}([x,y], \alpha(w),\alpha(z))+ 
  \epsilon(x, y+w)as_{T}([y,z], \alpha(w),\alpha(x))+  
  \epsilon(y, z+w)as_{T}([z,x], \alpha(w),\alpha(y))=0,$$
where the trilinear map $as_{T}: T \times T \times T \rightarrow T$ is the $Hom$-associator of $T$ defined by
$$as_{T}(x,y,z)=[[x,y],\alpha(z)]-[\alpha(x),[y,z]].$$
}
\end{df}


\subsection{$n$-ary $Hom$-algebras.}
In this subsection we give some definitions related to multiplicative $Hom$-algebras with $n$-ary multiplication.
In general  (see, \cite{cm15}), for the definition of $n$-ary $Hom$-algebras we must use  a family of homomorphisms $\{ \alpha_i\}_I,$
but $\alpha_i=\alpha$
and $\alpha([x_1, \ldots, x_n]) = [\alpha(x_1), \ldots, \alpha(x_n)]$ in the case of multiplicative $n$-ary $Hom$-algebras.

\begin{df}
A multiplicative $n$-ary $Hom$-Lie color algebra $(T, [\cdot, \ldots, \cdot], \epsilon, \alpha)$ is a $\mathbb{G}$-graded vector space  $T=\bigoplus_{g \in \mathbb{G}} T_g$ 
with an $n$-linear map $[\cdot,\ldots, \cdot]: T \times \ldots \times T \rightarrow T$ and 
an even linear map $\alpha: T \rightarrow T$ satisfying
$$[x_1, \ldots, x_i, x_{i+1}, \ldots, x_n]=- \epsilon( x_i, x_{i+1}) [x_1, \ldots, x_{i+1}, x_i, \ldots, x_n],$$
$$[\alpha (x_1), \ldots, \alpha (x_{n-1}),[y_1, \ldots, y_n]] =\sum_{i=1}^n \epsilon( X_{n-1}, Y_{i-1}) [\alpha (y_1), \ldots, \alpha (y_{i-1}),[x_1,\ldots, x_{n-1},y_i],\alpha (y_{i+1}),\ldots, \alpha (y_n)].$$
\end{df}


\begin{df}
A multiplicative  $n$-ary $Hom$-associative color algebra $(T, [\cdot, \ldots, \cdot],\epsilon,  \alpha)$ is a $\mathbb{G}$-graded vector space  $T=\bigoplus_{g \in \mathbb{G}} T_g$ 
with an $n$-linear map $[\cdot,\ldots, \cdot]: T \times \ldots \times T \rightarrow T$ and 
an   even linear map $\alpha: T \rightarrow T$ satisfying
$$[ \alpha(x_1), \ldots, \alpha (x_i), [x_{i+1}, \ldots, x_{i+n}], \alpha (x_{i+n+1}), \ldots, \alpha (x_{2n-1})]=
[ \alpha (x_1), \ldots, \alpha (x_{n-1}), [x_{n}, \ldots, x_{2n-1}]].$$
\end{df}

Now we define the notion of $n$-ary multiplicative $Hom$-$\Omega$ color algebra for arbitrary family of polynomial identities $\Omega.$

\begin{df}For a (maybe $n$-ary) multilinear polynomial $f(x_1, \ldots, x_n)$
we fix the order of indices $\{i_1, \ldots, i_n\}$ of one non-associative word $[x_{i_1}\ldots x_{i_n}]_{\beta}$
from the polynomial $f$.
Here, $f = \sum\limits_{\beta, \sigma \in S_n} \alpha_{ \beta, \sigma} [x_{\sigma(i_1)} \ldots x_{\sigma(i_n)}]_{\beta},$
where $\beta$ is an arrangement of brackets in the non-associative word.
For the  shift  $\mu_i: \{{j_1}, \ldots, {j_n} \} \mapsto \{j_1, \ldots j_{i+1}, j_i , \ldots,  j_n\}$
we define the element $\epsilon(x_{j_i},x_{j_{i+1}}).$
Now, for arbitrary non-associative word $ [x_{\sigma(i_1)} \ldots x_{\sigma(i_n)}]_{\beta}$ its order of indexes is a composition of suitable shifts $\mu_i$,
and for this word we set $\epsilon_{\sigma}$ defined as the product of corresponding $\epsilon(x_{j_i},x_{j_{i+1}}).$
Now, for the multilinear polynomial $f$, we define the color multilinear polynomial
$$f_{co}=\sum\limits_{\beta, \sigma \in S_n} \alpha_{\beta, \sigma} \epsilon_{\sigma} [x_{\sigma(i_1)} \ldots x_{\sigma(i_n)}]_{\beta}.$$
To construct a 
color multiplicative multilinear $Hom$-polynomial $f_{co}^{Hom}$ from a color  multilinear polynomial $f_{co}$ we use the algorithm from \cite{cm15} and we are changing all free 
letters $x_j$ to $\alpha(x_j)$ in all words of the 
color  multilinear polynomial $f_{co}$. 
\end{df}

\begin{df}Let $\Omega=\{ f_i\}$ be a family of $n$-ary multilinear polynimials. Then

(1) A color  $n$-ary  $\Omega$-algebra $L$ is a color $n$-ary algebra satisfying the family of color  multilinear polinomials $\Omega_{co}=\left\{ (f_i)_{co} \right\}.$

(2) A multiplicative $n$-ary  $Hom$-$\Omega$ color algebra $T$ is a color $n$-ary algebra satisfying the family of 
color multiplicative multilinear $Hom$-polynomials  $\Omega_{co}^{Hom}= \left\{ (f_i)^{Hom}_{co} \right\}$
for a homomorphism  $\alpha.$
\end{df}


\subsection{Linear maps on $Hom$-algebras.}
From now on, we denote a multiplicative $n$-ary $Hom$-$\Omega$ color algebra $(T,[\cdot, \ldots, \cdot], \epsilon, \alpha)$ by  $T.$

\begin{df}
 A $\mathbb{G}$-graded subspace $M \subseteq T$ is a color $Hom$-subalgebra of $T$ if $\alpha(M) \subseteq M$ and 
 $[M, \ldots, M] \subseteq  M.$ 
A $\mathbb{G}$-graded subspace $I \subseteq T$ is a color $Hom$-ideal of $T$ if $\alpha(I) \subseteq I$ and 
$[T, \ldots, I, \ldots, T] \subseteq I.$ 
\end{df}

\begin{df}
The vector subspace $S$ of ${\rm End}(T)$ is defined by $S =\{u \in {\rm End}(T) \ | \ u\alpha =\alpha u \}$ with $\tilde{\alpha}: S \rightarrow S$ given by $\tilde{\alpha}(u)=\alpha u$.
\end{df}

Notice that $(S,[\cdot,\cdot],\epsilon,\tilde{\alpha})$, where $[\cdot,\cdot]$ is the color bracket in (\ref{colorcom}) and $\epsilon$ is a bicharacter, is a color $Hom$-Lie algebra over $\mathbb{F}$.
In what follows, we only consider linear maps from S. Let $\xi \in \mathbb{G}.$

\begin{df}
$D \in {\rm End_{\xi}(T)}$ is a homogeneous  $\alpha^k$-derivation of $T$, with $k \in \mathbb{N}$ if it satisfies
$$D([x_1, \ldots, x_n])=\sum\limits_{i=1}^n \epsilon(D, X_{i-1})[\alpha^k(x_1), \ldots, \alpha^k(x_{i-1}), D(x_i), \alpha^k(x_{i+1}), \ldots, \alpha^k(x_n)].$$
\end{df}

We denote the set of all $\alpha^{k}$-derivations of $T$ by ${\rm Der}_{\alpha^k}(T)$. Notice that ${\rm Der}(T):=\bigoplus_{k \geq 0} {\rm Der}_{\alpha^k}(T)$ equipped with the color bracket defined by (\ref{colorcom}) and the even map $\tilde{\alpha}: {\rm Der}(T) \rightarrow {\rm Der}(T)$ defined by $\tilde{\alpha}(D)=D\alpha$,
is a color $Hom$-subalgebra of $S$, which is called the derivation algebra of $T$.

\begin{df}
 $D \in {\rm End_{\xi}(T)}$ is a homogeneous generalized $\alpha^k$-derivation of degree $\xi$ of $T$ if there exist endomorphisms $D^{(i)} \in {\rm End}_{\xi}(T)$ satisfying
$$
[D(x_1), \alpha^k(x_2), \ldots, \alpha^k(x_n)] + \sum\limits_{i=2}^{n} \epsilon(D, X_{i-1}) [\alpha^k(x_1), \ldots, \alpha^k(x_{i-1}), D^{(i-1)}(x_i), \alpha^k(x_{i+1}), \ldots, \alpha^k(x_n)]
$$
$$=D^{(n)}([x_1, \dots, x_n]).$$
\end{df}

\begin{df}
$D \in {\rm End_{\xi}(T)}$ is a homogeneous $\alpha^k$-quasiderivation of degree $\xi$ of $T$ if there exists $D' \in {\rm End}_{\xi}(T)$ satisfying
$$
 \sum\limits_{i=1}^{n} \epsilon(D, X_{i-1}) [\alpha^k(x_1), \ldots, \alpha^k(x_{i-1}), D(x_i), \alpha^k(x_{i+1}), \ldots, \alpha^k(x_n)]=D'([x_1, \dots, x_n]).$$
\end{df}

Let ${\rm GDer}_{\alpha^k}(T)$ and ${\rm QDer}_{\alpha^k}(T)$ be the sets of homogeneous generalized $\alpha^k$-derivations and of homogeneous $\alpha^k$-quasiderivations, respectively. We denote

\begin{center}
${\rm GDer}(T):= \bigoplus_{k \geq 0} {\rm GDer}_{\alpha^k}(T), \
  {\rm QDer}(T):= \bigoplus_{k \geq 0} {\rm QDer}_{\alpha^k}(T).$
\end{center}

\begin{df}
$D \in {\rm End}_{\xi}(T)$  is a homogeneous $\alpha^k$-centroid of $T$, where $k \in \mathbb{N}$ if it satisfies
$$D([x_1, \dots, x_n])=
\epsilon(D, X_{i-1}) [\alpha^k(x_1),\ldots,\alpha^k(x_{i-1}), D(x_i), \alpha^k(x_{i+1}), \ldots, \alpha^k(x_n)]$$
for all $i.$
\end{df}

\begin{df}
$D \in {\rm End}_{\xi}(T)$ is a homogeneous $\alpha^k$-quasicentroid of degree $\xi$ of $T$, where $k \in \mathbb{N}$ if it satisfies
$$
[D(x_1), \alpha^k(x_2), \ldots, \alpha^k(x_n)] =
\epsilon(D, X_{i-1}) [\alpha^k(x_1),\ldots,\alpha^k(x_{i-1}), D(x_i), \alpha^k(x_{i+1}), \ldots, \alpha^k(x_n)]
$$
for all $i.$
\end{df}

Denote the sets of homogeneous $\alpha^k$-centroids and homogeneous $\alpha^k$-quasicentroids by $C_{\alpha^k}(T)$ and $QC_{\alpha^k}(T)$, respectively. Denote also

\begin{center}
$C(T):= \bigoplus_{k \geq 0} C_{\alpha^k}(T), \ QC(T):= \bigoplus_{k \geq 0} QC_{\alpha^k}(T).$
\end{center}

\begin{df}$D \in {\rm End_{\xi}(T)}$ is a homogeneous  $\alpha^k$-center derivation of degree $\xi$ of $T$ if it satisfies

$$
 D([x_1, \ldots, x_n])=[\alpha^k(x_1), \ldots, \alpha^k(x_{i-1}), D(x_i), \alpha^k(x_{i+1}), \ldots, \alpha^k(x_n)]=0.
$$

\end{df}

Denote the set of homogeneous $\alpha^k$-center derivations by ${\rm ZDer}_{\alpha^k}(T)$. We also denote

\begin{center}
${\rm ZDer}(T):= \bigoplus_{k \geq 0} {\rm ZDer}_{\alpha^k}(T).$
\end{center}

Notice that we have the following chain of inclusions:
$${\rm ZDer}(T) \subseteq {\rm Der}(T) \subseteq {\rm QDer}(T) \subseteq {\rm GDer} (T) \subseteq {\rm End}(T).$$


\section{Generalized derivation algebras and their color $Hom$-subalgebras}

  In this section we  present some basic properties
of generalized derivations, quasiderivations and center derivations of a multiplicative  $n$-ary $Hom$-$\Omega$ color algebra.
\vspace{0,3cm}

\begin{Lem}
Let $T$ be a multiplicative  $n$-ary $Hom$-$\Omega$ color algebra. Then the following statements hold:

$(1)$\quad ${\rm GDer}(T),{\rm QDer}(T)$ and ${\rm C}(T)$ are color
Hom-subalgebras of $S$;

$(2)$\quad ${\rm ZDer}(T)$ is a color $Hom$-ideal of ${\rm Der}(T)$.
\end{Lem}

\rm {\it Proof.}\quad  $(1)$ Let $D_{\xi} \in{\rm
GDer}_{\alpha^k}(T),D_{\eta}\in{\rm
GDer}_{\alpha^s}(T)$, where $k, s \in \mathbb{N}$. For arbitrary $x_1, \ldots, x_n \in T$ we have

$[\tilde{\alpha}(D_{\xi})(x_1),\alpha^{k+1}(x_2),\ldots, \alpha^{k+1}(x_n)]
=[(D_{\xi}\alpha)(x_1),\alpha^{k+1}(x_2),\ldots, \alpha^{k+1}(x_n)]
=\alpha([D_{\xi}(x_1),\alpha^{k}(x_2),\ldots, \alpha^{k}(x_n)])$

$=\alpha(D_{\xi}^{(n)}([x_1, x_2,\ldots,x_n])-\sum\limits_{i=2}^n \epsilon(D_{\xi},X_{i-1})[\alpha^k(x_1),\ldots, \alpha^k(x_{i-1}),D_{\xi}^{(i-1)}(x_i),\alpha^k(x_{i+1}),\ldots, \alpha^k(x_n)])$

$=\tilde{\alpha}(D_{\xi}^{(n)})([x_1,x_2,\ldots,x_n])-\sum\limits_{i=2}^n \epsilon(D_{\xi},X_{i-1})[\alpha^{k+1}(x_1),\ldots, \alpha^{k+1}(x_{i-1}),\tilde{\alpha}(D_{\xi}^{(i-1)})(x_i),\alpha^{k+1}(x_{i+1}),\ldots, \alpha^{k+1}(x_n)]$.

For any $i \in \mathbb{N}$, $\tilde{\alpha}(D_{\xi}^{(i)})$ belongs to End$_{\xi}(T)$, thus $\tilde{\alpha}(D_{\xi}) \in $ GDer$_{\alpha^{k+1}}(T)$ and is obviously of degree $\xi$. For arbitrary $x_1, \ldots, x_n \in T$ we obtain

$[D_{\xi}D_{\eta}(x_1), \alpha^{k+s}(x_2), \ldots, \alpha^{k+s}(x_n)]$

$=D_{\xi}^{(n)}([D_{\eta}(x_1), \alpha^s(x_2), \ldots, \alpha^s(x_n)])$ 

$-\sum\limits_{i=2}^n \epsilon(D_{\xi},D_{\eta}+X_{i-1}) [\alpha^k(D_{\eta}(x_1)), \alpha^{k+s}(x_2),\ldots, \alpha^{k+s}(x_{i-1}),D_{\xi}^{(i-1)}(\alpha^s(x_i)),\alpha^{k+s}(x_{i+1}), \ldots, \alpha^{k+s}(x_n)]$

$=D_{\xi}^{(n)}D_{\eta}^{(n)}([x_1,\ldots, x_n])$

$-\sum\limits_{j=2}^n \epsilon(D_{\eta},X_{j-1}) D_{\xi}^{(n)}[\alpha^s(x_1),\ldots, \alpha^s(x_{j-1}), D_{\eta}^{(j-1)}(x_{j}), \alpha^s(x_{j+1}), \ldots, \alpha^s(x_n)]$

$-\sum\limits_{i=2}^n \epsilon(D_{\xi},D_{\eta}+X_{i-1}) D_{\eta}^{(n)}[\alpha^k(x_1),\ldots, \alpha^k(x_{i-1}), D_{\xi}^{(i-1)}(x_{i}), \alpha^k(x_{i+1}), \ldots, \alpha^k(x_n)]$

$+\sum\limits_{i=2, j=2, j<i}^n \epsilon(D_{\xi},D_{\eta}+X_{i-1})\epsilon(D_{\eta}, X_{j-1}) [\alpha^{k+s}(x_1),\ldots, D_{\eta}^{(j-1)}(\alpha^k(x_{j})), \ldots, D_{\xi}^{(i-1)}(\alpha^s(x_i)), \ldots, \alpha^{k+s}(x_{n})]$

$+\sum\limits_{i=2, j=2, i<j}^n \epsilon(D_{\xi},X_{i-1})\epsilon(D_{\eta}, X_{j-1}) [\alpha^{k+s}(x_1),\ldots, D_{\xi}^{(i-1)}(\alpha^s(x_{i})), \ldots, D_{\eta}^{(j-1)}(\alpha^k(x_{j})), \ldots, \alpha^{k+s}(x_n)]$

$+\sum\limits_{i=2}^n \epsilon(D_{\xi},D_{\eta} + X_{i-1}) \epsilon(D_{\eta}, X_{i-1}) [\alpha^{k+s}(x_1), \ldots, \alpha^{k+s}(x_{i-1}), D_{\eta}^{(i-1)}D_{\xi}^{(i-1)}(x_i), \alpha^{k+s}(x_{i+1}), \ldots, \alpha^{k+s}(x_n)]$.

Thus, for arbitrary $x_1, \ldots, x_n \in T$, we have

$[[D_{\xi},D_{\eta}](x_1),\alpha^{k+s}(x_2), \ldots, \alpha^{k+s}(x_n)]
=[(D_{\xi}D_{\eta}-\epsilon(D_{\xi},D_{\eta})D_{\eta}D_{\xi})(x_1), \alpha^{k+s}(x_2), \ldots, \alpha^{k+s}(x_n)]$

$=[D_{\xi}^{(n)}, D_{\eta}^{(n)}]([x_1, \ldots, x_n])-\sum\limits_{i=2}^{n} \epsilon(D_{\xi}+D_{\eta}, X_{i-1}) [\alpha^{k+s}(x_1), \ldots, [D_{\xi}^{(i-1)},D_{\eta}^{(i-1)}](x_i),\ldots, \alpha^{k+s}(x_n)]$.

For all $2 \leq i \leq n$, $[D_{\xi}^{(i-1)}, D_{\eta}^{(i-1)}] \in$ ${\rm End_{\xi+\eta}}(T)$. Therefore $[D_{\xi}^{(n)}, D_{\eta}^{(n)}] \in$ ${\rm GDer_{\alpha^{k+s}}}(T)$ and ${\rm GDer}(T)$ is a color $Hom$-subalgebra of $S$.

Similarly, ${\rm QDer}(T)$ is a color $Hom$-subalgebra of $S$.

Let $D_{\xi} \in C_{\alpha^k}(T), D_{\eta} \in C_{\alpha^s}(T)$. 
For arbitrary $x_1, \ldots, x_n \in T$ we have

$\tilde{\alpha}(D_{\xi})([x_1, \ldots, x_n])=\alpha D_{\xi}([x_1, \ldots, x_n])
= \alpha([D_{\xi}(x_1), \alpha^k(x_2), \ldots, \alpha^k(x_n)])
= [\tilde{\alpha}({D_{\xi}})(x_1),\alpha^{k+1}(x_2), \ldots, \alpha^{k+1}(x_n)]$

and

$[\tilde{\alpha}({D_{\xi}})(x_1),\alpha^{k+1}(x_2),\ldots, \alpha^{k+1}(x_n)]
= \alpha([D_{\xi}(x_1), \alpha^k(x_2), \ldots, \alpha^k(x_n)])$

$= \epsilon(D_{\xi}, X_{i-1}) \alpha([\alpha^{k}(x_1), \ldots, \alpha^k(x_{i-1}), D_{\xi}(x_i), \alpha^k(x_{i+1}), \ldots, \alpha^k(x_n)])$

$= \epsilon(D_{\xi}, X_{i-1}) [\alpha^{k+1}(x_1), \ldots, \alpha^{k+1}(x_{i-1}), \tilde{\alpha}(D_{\xi})(x_i), \alpha^{k+1}(x_{i+1}), \ldots, \alpha^{k+1}(x_n)].$

Thus, $\tilde{\alpha}(D_{\xi}) \in C_{\alpha^{k+1}}(T)$.

Notice that

$[[D_{\xi}, D_{\eta}](x_1),\alpha^{k+s}(x_2),\ldots, \alpha^{k+s}(x_n)]$

$=[D_{\xi}D_{\eta}(x_1), \alpha^{k+s}(x_2), \ldots, \alpha^{k+s}(x_n)]-\epsilon(D_{\xi}, D_{\eta})[D_{\eta}D_{\xi}(x_1),\alpha^{k+s}(x_2),\ldots, \alpha^{k+s}(x_n)]$

$=D_{\xi}D_{\eta}([x_1,x_2, \ldots, x_n])- \epsilon(D_{\xi},D_{\eta}) D_{\eta}D_{\xi} ([x_1, x_2, \ldots, x_n])
=[D_{\xi}, D_{\eta}]([x_1, x_2, \ldots, x_n]).$

Similarly, we have

$\epsilon(D_{\xi}+D_{\eta},X_{i-1})[\alpha^{k+s}(x_1), \ldots, [D_{\xi},D_{\eta}](x_i), \ldots, \alpha^{k+s}(x_n)]=[D_{\xi},D_{\eta}]([x_1,x_2,\ldots, x_n])$.

Thus, $[D_{\xi},D_{\eta}] \in C_{\alpha^{k+s}}(T)$ of degree $\xi + \eta$ and $C(T)$ is a color $Hom$-subalgebra of $S$.

(2) Let $D_{\xi} \in ZDer_{\alpha^k}(T), D_{\eta} \in Der_{\alpha^s}(T)$. For arbitrary $x_1, \ldots, x_n \in T$ we have

$[\tilde{\alpha}(D_{\xi})(x_1), \alpha^{k+1}(x_2), \ldots, \alpha^{k+1}(x_n)]
=[(\alpha D_{\xi})(x_1), \alpha^{k+1}(x_2),\ldots, \alpha^{k+1}(x_n)]$

$=\alpha([D_{\xi}(x_1), \alpha^k(x_2), \ldots, \alpha^k(x_n)])
=\alpha D_{\xi} ([x_1, x_2, \ldots, x_n])=\tilde{\alpha}(D_{\xi})([x_1, x_2, \ldots, x_n])=0.$

So, $\tilde{\alpha}(D_{\xi}) \in ZDer_{\alpha^{k+1}}(T)$.

Observe that

$[D_{\xi},D_{\eta}]([x_1, x_2, \ldots, x_n])
= D_{\xi}D_{\eta} ([x_1, x_2, \ldots, x_n])-\epsilon(D_{\xi},D_{\eta}) D_{\eta}D_{\xi}([x_1,x_2, \ldots, x_n])$

$=D_{\xi}([D_{\eta}(x_1), \alpha^s(x_2), \ldots, \alpha^s(x_n)]) +\sum\limits_{i=2}^n \epsilon(D_{\eta}, X_{i-1})D_{\xi}( [\alpha^s(x_1),\ldots, \alpha^s(x_{i-1}), D_{\eta}^{(i-1)}(x_i), \alpha^s(x_{i+1}), \ldots, \alpha^s(x_n)])=0$

and

$[[D_{\xi},D_{\eta}](x_1), \alpha^{k+s}(x_2), \ldots, \alpha^{k+s}(x_n)]$

$=[D_{\xi}D_{\eta}(x_1), \alpha^{k+s}(x_2), \ldots, \alpha^{k+s}(x_n)]-\epsilon(D_\xi, D_{\eta}) [D_{\eta}D_{\xi}(x_1), \alpha^{k+s}(x_2), \ldots, \alpha^{k+s}(x_n)]$

$=-\epsilon(D_\xi, D_{\eta}) [D_{\eta}(D_{\xi}(x_1)), \alpha^{k+s}(x_2), \ldots, \alpha^{k+s}(x_n)]
=-\epsilon(D_\xi, D_{\eta}) D_{\eta}([D_{\xi}(x_1), \alpha^{k}(x_2), \ldots, \alpha^{k}(x_n)])$

$+\sum\limits_{i=2}^{n} \epsilon(D_\xi, D_{\eta})\epsilon(D_{\eta}, D_{\xi}+X_{i-1}) [\alpha^s(D_{\xi}(x_1)),\alpha^{k+s}(x_2), \ldots, D_{\eta}^{(i-1)}(\alpha^k(x_i)), \ldots, \alpha^{k+s}(x_n)]=0.$ 

Hence $[D_{\xi}, D_{\eta}] \in ZDer_{\alpha^{k+s}}(T)$ and is of degree $\xi + \eta$. Therefore, $ZDer(T)$ is a color $Hom$-ideal of ${\rm Der} (T)$.
\vspace{0,3cm}

\begin{Lem}\label{22} Let
$T$ be a multiplicative  $n$-ary $Hom$-$\Omega$ color algebra. Then

    $(1)$ \quad $[{\rm Der}(T),{\rm C}(T)]\subseteq {\rm C}(T);$

    $(2)$ \quad $[{\rm QDer}(T),{\rm QC}(T)]\subseteq {\rm QC}(T);$

    $(3)$ \quad ${\rm C}(T)\cdot{\rm Der}(T)\subseteq {\rm Der}(T);$

    $(4)$ \quad ${\rm C}(T)\subseteq {\rm QDer}(T);$

    $(5)$ \quad $[{\rm QC}(T),{\rm QC}(T)]\subseteq {\rm QDer}(T);$

    $(6)$ \quad ${\rm QDer}(T)+{\rm QC}(T)\subseteq {\rm GDer}(T).$
\end{Lem}

 \rm {\it Proof.}\quad $(1)$ Let $D_{\xi} \in {\rm Der}_{\alpha^k}(T), D_{\eta} \in {\rm C}_{\alpha^s}(T)$. For arbitrary $x_1, \ldots, x_n \in T$ we have
 
$[D_{\xi}D_{\eta}(x_1), \alpha^{k+s}(x_2), \ldots, \alpha^{k+s}(x_n)]$

$=D_{\xi}([D_{\eta}(x_1), \alpha^s(x_2), \ldots, \alpha^s(x_n)])
-\sum\limits_{i=2}^n \epsilon(D_{\xi}, D_{\eta} + X_{i-1})[D_{\eta}(\alpha^k(x_1)), \alpha^{k+s}(x_2), \ldots, D_{\xi}(\alpha^s(x_i)), \ldots, \alpha^{k+s}(x_n)]$

$=D_{\xi}D_{\eta}([x_1, x_2, \ldots, x_n])
-\sum\limits_{i=2}^n \epsilon(D_{\xi}, D_{\eta}) \epsilon(D_{\xi}+ D_{\eta},X_{i-1}) [\alpha^{k+s}(x_1), \ldots, D_{\eta}D_{\xi}(x_i), \ldots, \alpha^{k+s}(x_n)]$

and

$[D_{\eta}D_{\xi}(x_1), \alpha^{k+s}(x_2), \ldots, \alpha^{k+s}(x_n)]$

$=D_{\eta}([D_{\xi}(x_1), \alpha^k(x_2), \ldots, \alpha^k(x_n)])$

$=D_{\eta}D_{\xi}([x_1, x_2, \ldots, x_n])
-\sum\limits_{i=2}^n \epsilon(D_{\xi}, X_{i-1})D_{\eta}([\alpha^k(x_1),  \ldots, D_{\xi}(x_i), \ldots, \alpha^{k}(x_n)])$

$=D_{\eta}D_{\xi}([x_1, x_2, \ldots, x_n])
-\sum\limits_{i=2}^n \epsilon(D_{\xi}+D_{\eta}, X_{i-1})[\alpha^{k+s}(x_1),  \ldots, D_{\eta}D_{\xi}(x_i), \ldots, \alpha^{k+s}(x_n)]$.

Hence

$[[D_{\xi}, D_{\eta}](x_1), \alpha^{k+s}(x_2), \ldots, \alpha^{k+s}(x_n)]$

$=[D_{\xi}D_{\eta}(x_1), \alpha^{k+s}(x_2), \ldots, \alpha^{k+s}(x_n)]-\epsilon(D_{\xi},D_{\eta}) [D_{\eta}D_{\xi}(x_1), \alpha^{k+s}(x_2), \ldots, \alpha^{k+s}(x_n)]$

$= [D_{\xi}, D_{\eta}]([x_1, \ldots, x_n])$.

Similarly, 

$[[D_{\xi}, D_{\eta}](x_1), \alpha^{k+s}(x_2), \ldots, \alpha^{k+s}(x_n)]= \epsilon(D_{\xi}+D_{\eta}, X_{i-1})[\alpha^{k+s}(x_1), \ldots, [D_{\xi},D_{\eta}](x_i), \ldots, \alpha^{k+s}(x_n)]$. 

Thus, $[D_{\xi},D_{\eta}] \in {\rm C}_{\alpha^{k+s}}(T)$, is of degree $\xi + \eta$ and $[{\rm Der}(T), {\rm C}(T)] \subseteq {\rm C}(T)$.

$(2)$ Is similar to the proof of $(1)$.

$(3)$ Let $D_{\xi}\in{\rm C_{\alpha^k}}(T), D_{\eta}\in{\rm Der}_{\alpha^s}(T)$. For arbitrary
$x_1,\ldots, x_n \in T$ we have

$D_{\xi}D_{\eta}([x_1, \ldots, x_n])$

$=D_{\xi}(\sum\limits_{i=1}^n \epsilon (D_{\eta}, X_{i-1}) [\alpha^{s}(x_1), \ldots, D_{\eta}(x_i), \ldots, \alpha^{s}(x_n)])$

$=\sum\limits_{i=1}^n \epsilon (D_{\xi}+D_{\eta}, X_{i-1}) [\alpha^{k+s}(x_1), \ldots, D_{\xi}D_{\eta}(x_i), \ldots, \alpha^{k+s}(x_n)]$.

Therefore $D_{\xi}D_{\eta}\in {\rm Der}_{\alpha^{k+s}}(T)$ and is of degree $\xi + \eta$. Hence ${\rm C}(T)\cdot{\rm Der}(T)\subseteq {\rm Der}(T)$.

$(4)$ Let $D_{\xi}\in{\rm C}_{\alpha^k}(T)$. For arbitrary $x_1,\ldots, x_n \in T$ we have

$D_{\xi}([x_1, \ldots, x_n])
= \epsilon(D_{\xi},X_{i-1})[\alpha^{k}(x_1), \ldots, D_{\xi}(x_i),\ldots, \alpha^{k}(x_n)].$

Hence $\sum\limits_{i=1}^n  \epsilon(D_{\xi},X_{i-1})[\alpha^{k}(x_1), \ldots, D_{\xi}(x_i),\ldots, \alpha^{k}(x_n)]=nD_{\xi}[x_1,\ldots, x_n].$

Therefore $D_{\xi}\in {\rm QDer}_{\alpha^k}(T)$ since $D^{'}=nD_{\xi}\in {\rm
C}_{\alpha^k}(T)$.

$(5)$ Let $D_{\xi} \in {\rm QC}_{\alpha^k} (T), D_{\eta} \in {\rm QC}_{\alpha^s}(T).$ For arbitrary $x_1, \ldots, x_n \in T$
we have

$[\alpha^{k+s}(x_1), \ldots, [D_{\xi},D_{\eta}](x_i),\ldots, \alpha^{k+s}(x_n)]$

$=\epsilon(X_{i-1},D_{\xi})\epsilon(X_{i-1}-X_1,D_{\eta})[D_{\xi}(\alpha^s(x_1)),D_{\eta}(\alpha^k(x_2)), \ldots, \alpha^{k+s}(x_i),\ldots, \alpha^{k+s}(x_n)]$

$-\epsilon(D_{\xi},D_{\eta})\epsilon(X_{i-1}-X_1,D_{\eta})\epsilon(X_{i-1}+D_{\eta},D_{\xi}) [D_{\xi}(\alpha^s(x_1)),D_{\eta}(\alpha^k(x_2)), \ldots, \alpha^{k+s}(x_i), \ldots, \alpha^{k+s}(x_n)]=0.$

Hence

$\sum\limits_{i=1}^n \epsilon(D_{\xi}+D_{\eta}, X_{i-1}) [\alpha^{k+s}(x_1), \ldots, [D_{\xi},D_{\eta}](x_i), \ldots, \alpha^{k+s}(x_n)]=0.$

therefore $[D_{\xi},D_{\eta}]\in {\rm QDer}_{\alpha^{k+s}}(T)$ and is of degree $\xi+\eta$.

$(6)$\quad Is obvious.\hfill$\Box$

\begin{Lem}If $T$ is a multiplicative  $n$-ary $Hom$-$\Omega$ color algebra,
then ${\rm QC}(T)+[{\rm QC}(T),{\rm QC}(T)]$ is a color $Hom$-subalgebra of
${\rm GDer}(T)$.
\end{Lem}

{\it Proof.} \quad By Lemma \ref{22}~$(5)$ and $(6)$, 
we have

${\rm QC}(T)+[{\rm QC}(T),{\rm QC}(T)]\subseteq {\rm GDer}(T)$ and

$ [{\rm QC}(T)+[{\rm QC}(T),{\rm QC}(T)],{\rm QC}(T)+[{\rm QC}(T),{\rm QC}(T)]]$

$\subseteq [{\rm QC}(T)+{\rm QDer}(T),{\rm QC}(T)+[{\rm QC}(T),{\rm QC}(T)]]$

$\subseteq [{\rm QC}(T),{\rm QC}(T)]+[{\rm QC}(T),[{\rm QC}(T),{\rm
QC}(T)]]+[{\rm QDer}(T),{\rm QC}(T)] + [{\rm QDer}(T),[{\rm
QC}(T),{\rm QC}(T)]]$ 

Using the color $Hom$-Jacobi identity and $(2)$ of the previous lemma, it is easy to verify that 

$[{\rm
QDer}(T),[{\rm QC}(T),{\rm QC}(T)]]\subseteq [{\rm QC}(T),{\rm
QC}(T)]$. 

Thus, 

$[{\rm
QC}(T)+[{\rm QC}(T),{\rm QC}(T)],{\rm QC}(T)+[{\rm QC}(T),{\rm
QC}(T)]]\subseteq {\rm QC}(T)+[{\rm QC}(T),{\rm QC}(T)]$.
\hfill$\Box$

\begin{df}
Let $T$ is a multiplicative $n$-ary algebra.
 The annihilator   of $T$ is defined by
$${\rm Ann}(T) =  \{ x \in T | \ \  [ x_1, \ldots, x_{i-1}, x, x_{i+1}, \ldots, x_{n}]=0 \text{ {\rm for all }} i\}.$$
\end{df}

\begin{Lem}Let $\alpha$ be a surjective map. If $T$ is a multiplicative  $n$-ary $Hom$-$\Omega$ color algebra,
then
$[{\rm C}(T),{\rm QC}(T)]\subseteq {\rm End}(T,{\rm Ann}(T))$.
 Moreover, if ${\rm Ann}(T)=\{0\},$ then $[{\rm C}(T),{\rm QC}(T)]=\{0\}.$
\end{Lem}

{\it Proof.}\quad Assume that
$D_{\xi}\in {\rm C}_{\alpha^k}(T), D_{\eta}\in {\rm QC}_{\alpha^s}(T)$ and $x \in T$. As $\alpha$ is surjective, for any $y_i' \in T, 1 \leq i \leq n$, there exists $y_i \in T$ such that $y_i' = \alpha^{k+s} (y_i)$.
We have


$\epsilon(Y_{i-1},D_{\xi}+D_{\eta})[\alpha^{k+s}(y_1), \ldots, [D_{\xi},D_{\eta}](x), \ldots, \alpha^{k+s}(y_n)]=$

$\epsilon(Y_{i-1},D_{\xi}+D_{\eta})[\alpha^{k+s}(y_1), \ldots, (D_{\xi}D_{\eta} - \epsilon(D_{\xi},D_{\eta}) D_{\eta}D_{\xi})(x), \ldots, \alpha^{k+s}(y_n)]=$



$D_{\xi}([D_{\eta}(y_1), \alpha^{s}(y_2), \ldots, \alpha^{s}(x), \ldots,\alpha^{s}(y_n)]-
[D_{\eta}(y_1),\alpha^s(y_2), \ldots, \alpha^{s}(x), \ldots, \alpha^{s}(y_n)])=0$.

Hence
$[D_{\xi},D_{\eta}](x)\in {\rm Ann}(T)$ and $[D_{\xi},D_{\eta}]\in {\rm
End}(T,{\rm Ann}(T))$ as desired. Furthermore, if ${\rm Ann}(T)=\{0\},$
it is clear that $[{\rm C}(T),{\rm QC}(T)]=\{0\}.$
\hfill$\Box$  

\begin{Lem}\label{27}Let $\mathbb{F}$ be of characteristic $\neq 2$ and $(S, \bullet, \tilde{\alpha})$ be a color $Hom$-algebra
with multiplication 

\begin{center}
$D_{\xi} \bullet D_{\eta} = \frac 1 2 (D_{\xi}D_{\eta} +\epsilon(D_{\xi},D_{\eta})D_{\eta}D_{\xi})$
\end{center}
where $D_{\xi}, D_{\eta}$ are $\alpha$-derivations of $S$. Then

$(1)$ \quad $(S, \bullet, \alpha)$ is a color $Hom$-Jordan algebra;

$(2)$ \quad $({\rm QC}(T), \bullet, \alpha)$ is a color $Hom$-Jordan algebra. 

\end{Lem}

{\it Proof.}\quad $(1)$ Let $D_{\xi}, D_{\eta} \in S$. We have

$D_{\xi} \bullet D_{\eta}$

$= \frac 1 2 (D_{\xi}D_{\eta}+\epsilon(D_{\xi},D_{\eta})D_{\eta}D_{\xi})$

$= \frac 1 2 \epsilon(D_{\xi},D_{\eta}) (D_{\eta}D_{\xi}+\epsilon(D_{\eta},D_{\xi})D_{\xi}D_{\eta})$

$= \epsilon(D_{\xi}, D_{\eta}) D_{\eta} \bullet D_{\xi}$.

Consider now the second $Hom$-Jordan identity:

$((D_{\xi} \bullet D_{\eta}) \bullet \alpha(D_{\theta})) \bullet \alpha^2(D_{\gamma})$

$= \frac 1 2 ((D_{\xi}D_{\eta} + \epsilon(D_{\xi}, D_{\eta}) D_{\eta} D_{\xi}) \bullet \alpha(D_{\theta}))\bullet \alpha^2(D_{\gamma})$  

$=\frac 1 4 ((D_{\xi}D_{\eta} + \epsilon(D_{\xi}, D_{\eta}) D_{\eta} D_{\xi}) \alpha(D_{\theta}) + \epsilon(D_{\xi}+D_{\eta}, D_{\theta}) \alpha(D_{\theta})(D_{\xi}D_{\eta} + \epsilon(D_{\xi}, D_{\eta}) D_{\eta} D_{\xi})) \bullet \alpha^2(D_{\gamma})$

$=\frac 1 8 (D_{\xi}D_{\eta} \alpha(D_{\theta}) \alpha^2(D_{\gamma}) + \epsilon(D_{\xi}, D_{\eta}) D_{\eta} D_{\xi} \alpha(D_{\theta}) \alpha^2(D_{\gamma}) + \epsilon(D_{\xi}+D_{\eta},D_{\theta})\alpha(D_{\theta}) D_{\xi}D_{\eta} \alpha^2(D_{\gamma})$

$ + \epsilon(D_{\xi}+D_{\eta}, D_{\theta}) \epsilon(D_{\xi},D_{\eta}) \alpha(D_{\theta}) D_{\eta} D_{\xi} \alpha^2(D_{\gamma}) + \epsilon(D_{\xi}+D_{\eta}+D_{\theta}, D_{\gamma}) \alpha^2(D_{\gamma}) D_{\xi} D_{\eta} \alpha(D_{\theta})$

$ + \epsilon(D_{\xi},D_{\eta}) \epsilon(D_{\xi}+D_{\eta}+D_{\theta}, D_{\gamma}) \alpha^2(D_{\gamma}) D_{\eta} D_{\xi} \alpha(D_{\theta}) + \epsilon(D_{\xi}+D_{\eta},D_{\theta}) \epsilon(D_{\xi}+D_{\eta}+D_{\theta}, D_{\gamma}) \alpha^2(D_{\gamma}) \alpha(D_{\theta}) D_{\xi} D_{\eta}$

$+\epsilon(D_{\xi}+D_{\eta},D_{\theta}) \epsilon(D_{\xi}, D_{\eta}) \epsilon(D_{\xi}+D_{\eta}+D_{\theta}, D_{\gamma}) \alpha^2(D_{\gamma}) \alpha(D_{\theta}) D_{\eta} D_{\xi})$.

On the other hand, we have

$\alpha(D_{\xi} \bullet D_{\eta}) \bullet (\alpha (D_{\theta}) \bullet \alpha (D_{\gamma}))$

$= \frac 1 4 \alpha(D_{\xi}D_{\eta} + \epsilon(D_{\xi}, D_{\eta}) D_{\eta} D_{\xi}) \bullet (\alpha(D_{\theta})\alpha(D_{\gamma}) + \epsilon(D_{\theta}, D_{\gamma}) \alpha(D_{\gamma}) \alpha(D_{\theta}))$

$= \frac 1 8 (\alpha(D_{\xi}D_{\eta}) \alpha(D_{\theta})\alpha(D_{\gamma}) + \epsilon(D_{\xi}+D_{\eta}, D_{\theta}+D_{\gamma}) \alpha(D_{\theta}) \alpha(D_{\gamma}) \alpha(D_{\xi}D_{\eta})$

$ + \epsilon(D_{\theta},D_{\gamma}) \alpha(D_{\xi}D_{\eta})\alpha(D_{\gamma}) \alpha(D_{\theta}) + \epsilon(D_{\theta}, D_{\gamma}) \epsilon(D_{\xi}+D_{\eta}, D_{\theta}+D_{\gamma}) \alpha(D_{\gamma}) \alpha(D_{\theta}) \alpha(D_{\xi}D_{\eta})$

$+ \epsilon(D_{\xi}, D_{\eta}) \alpha(D_{\eta} D_{\xi}) \alpha(D_{\theta}) \alpha(D_{\gamma}) + \epsilon(D_{\xi},D_{\eta}) \epsilon(D_{\xi}+D_{\eta}, D_{\theta}+D_{\gamma}) \alpha(D_{\theta}) \alpha(D_{\gamma}) \alpha(D_{\eta} D_{\xi}) $

$+ \epsilon(D_{\xi}, D_{\eta}) \epsilon(D_{\theta}, D_{\gamma}) \alpha(D_{\eta}D_{\xi}) \alpha(D_{\gamma}) \alpha(D_{\theta})
+ \epsilon(D_{\xi},D_{\eta}) \epsilon(D_{\theta}, D_{\gamma}) \epsilon(D_{\xi}+D_{\eta}, D_{\gamma}+D_{\theta}) \alpha(D_{\gamma}) \alpha(D_{\theta}) \alpha(D_{\eta}D_{\xi}))$.

So, 

$\epsilon(D_{\gamma}, D_{\xi}+D_{\theta}) as_{\alpha}(D_{\xi} \bullet D_{\eta}, \alpha(D_{\theta}), \alpha(D_{\gamma}))$

$= \frac 1 8 \epsilon(D_{\gamma}, D_{\xi}+D_{\theta})(\epsilon(D_{\xi}+D_{\eta},D_{\theta})\alpha(D_{\theta})D_{\xi} D_{\eta} \alpha^2(D_{\gamma}) + \epsilon(D_{\xi}, D_{\eta}) \epsilon(D_{\xi}+D_{\eta}, D_{\theta}) \alpha(D_{\theta}) D_{\eta} D_{\xi} \alpha^2(D_{\gamma})$

$+\epsilon(D_{\xi}+D_{\eta}+D_{\theta}, D_{\gamma}) \alpha^2(D_{\gamma}) D_{\xi} D_{\eta} \alpha(D_{\theta}) + \epsilon(D_{\xi}, D_{\eta})\epsilon(D_{\xi}+D_{\eta}+D_{\theta}, D_{\gamma}) \alpha^2(D_{\gamma}) D_{\eta} D_{\xi} \alpha(D_{\theta})$

$-\epsilon(D_{\xi}+D_{\eta}, D_{\theta}+D_{\gamma}) \alpha(D_{\theta}) \alpha(D_{\gamma}) \alpha(D_{\xi}D_{\eta})-\epsilon(D_{\theta},D_{\gamma})\alpha(D_{\xi}D_{\eta}) \alpha(D_{\gamma}) \alpha(D_{\theta})$

$-\epsilon(D_{\xi},D_{\eta}) \epsilon(D_{\xi}+D_{\eta}, D_{\theta}+D_{\gamma}) \alpha(D_{\theta}) \alpha(D_{\gamma}) \alpha(D_{\eta}D_{\xi})-\epsilon(D_{\xi},D_{\eta})\epsilon(D_{\theta},D_{\gamma})\alpha(D_{\eta}D_{\xi}) \alpha(D_{\gamma}) \alpha(D_{\theta}))$.

Hence

$\epsilon(D_{\gamma}, D_{\xi}+D_{\theta}) as_{\alpha}(D_{\xi} \bullet D_{\eta}, \alpha(D_{\theta}), \alpha(D_{\gamma})) + \epsilon(D_{\xi}, D_{\eta}+D_{\theta}) as_{\alpha}(D_{\eta} \bullet D_{\gamma}, \alpha(D_{\theta}), \alpha(D_{\xi}))$

$ + \epsilon(D_{\eta}, D_{\gamma}+D_{\theta}) as_{\alpha}(D_{\gamma} \bullet D_{\xi}, \alpha(D_{\theta}), \alpha(D_{\eta}))$=0.

$(2)$ We only need to show that for arbitrary $D_{\xi},D_{\eta}\in{\rm QC}(T)$, $D_{\xi}\bullet
D_{\eta}\in{\rm QC}(T)$. We have

$[D_{\xi}\bullet D_{\eta}(x_1), \alpha^{k+s}(x_2), \ldots, \alpha^{k+s}(x_n)]$

$= \frac 1 2 [D_{\xi}D_{\eta}(x_1), \alpha^{k+s}(x_2), \ldots, \alpha^{k+s}(x_n)]+ \frac 1 2 \epsilon(D_{\xi},D_{\eta}) [D_{\eta}D_{\xi}(x_1), \alpha^{k+s}(x_2), \ldots, \alpha^{k+s}(x_n)]$

$=\frac 1 2 \epsilon(D_{\xi},D_{\eta}+X_{i-1})[D_{\eta}(\alpha^k(x_1)), \alpha^{k+s}(x_2), \ldots, D_{\xi}(x_i), \ldots, \alpha^{k+s}(x_n)]$

$+ \frac 1 2 \epsilon(D_{\eta},X_{i-1})[D_{\xi}(\alpha^s(x_1)), \alpha^{k+s}(x_2), \ldots, D_{\eta}(x_i), \ldots, \alpha^{k+s}(x_n)]$

$=\frac 1 2 \epsilon(D_{\xi},D_{\eta}+X_{i-1})\epsilon(D_{\eta},X_{i-1})[\alpha^{k+s}(x_1), \ldots, D_{\eta}D_{\xi}(x_i), \ldots, \alpha^{k+s}(x_n)]$

$+ \frac 1 2 \epsilon(D_{\xi} + D_{\eta},X_{i-1})[\alpha^{k+s}(x_1),\ldots, D_{\xi}D_{\eta}(x_i), \ldots, \alpha^{k+s}(x_n)]$

$=\epsilon(D_{\xi}+D_{\eta},X_{i-1})[ \alpha^{k+s}(x_1),\ldots, D_{\xi}\bullet D_{\eta}(x_i), \ldots,  \alpha^{k+s}(x_n)].$

Then $D_{\xi}\bullet D_{\eta}\in {\rm QC}(T)$ and ${\rm QC}(T)$ is a color $Hom$-Jordan algebra.
\hfill$\Box$

\begin{Th} Let $T$ be a multiplicative  $n$-ary $Hom$-$\Omega$ color algebra. 
Then the following statements hold:

$(1)$\quad ${\rm QC}(T)$ is a color $Hom$-Lie algebra with $[D_{\xi},D_{\eta}]= D_{\xi}D_{\eta}-\epsilon(D_{\xi},D_{\eta})D_{\eta}D_{\xi}$ if and only if
${\rm QC}(T)$ is a color $Hom$-associative algebra with respect to the usual composition of operators;

    $(2)$\quad If char $\mathbb {F}$ is not $2,$
$\alpha$ is a surjective map and ${\rm Ann}(T)=\{0\},$ then
${\rm QC}(T)$ is a color $Hom$-Lie algebra if and only if
$~[{\rm QC}(T),{\rm QC}(T)]=\{0\}.$
\end{Th}

\rm {\it Proof.}\quad $(1)$ $(\Leftarrow)$ For arbitrary $D_{\xi} \in
{\rm QC}_{\alpha^k}(T),D_{\eta}\in
{\rm QC}_{\alpha^s}(T)$, we have $D_{\xi}D_{\eta} \in {\rm QC}_{\alpha^{k+s}}(T)$ and $D_{\eta}D_{\xi}\in
{\rm QC}_{\alpha^{k+s}}(T)$. So, $[D_{\xi},D_{\eta}]=D_{\xi}D_{\eta}-\epsilon(D_{\xi},D_{\eta})D_{\eta}D_{\xi}\in {\rm
QC}_{\alpha^{k+s}}(T).$ Hence, ${\rm QC}(T)$ is a color $Hom$-Lie algebra.

$(\Rightarrow)$ Note that $D_{\xi}D_{\eta}=D_{\xi}\bullet D_{\eta}+\frac{[D_{\xi},D_{\eta}]}{2}$. 
By Lemma \ref{27}, we have
$D_{\xi}\bullet D_{\eta}\in {\rm QC}(T), [D_{\xi},D_{\eta}]\in {\rm QC}(T)$.
It follows that $D_{\xi}D_{\eta}\in {\rm QC}(T)$, as desired.

$(2)$ $(\Rightarrow)$\quad  Let $D_{\xi} \in {\rm QC}_{\alpha^k}(T), D_{\eta}\in {\rm QC}_{\alpha^s}(T)$.
Now we  have

$[[D_{\xi},D_{\eta}](x_1), \alpha^{k+s}(x_2), \ldots, \alpha^{k+s}(x_n)]$

$=\epsilon(D_{\xi},D_{\eta}+X_{i-1})\epsilon(D_{\eta},X_{i-1})[\alpha^{k+s}(x_1), \ldots, D_{\eta}D_{\xi}(x_i), \ldots, \alpha^{k+s}(x_n)]$

$- \epsilon(D_{\xi} + D_{\eta},X_{i-1})[\alpha^{k+s}(x_1),\ldots, D_{\xi}D_{\eta}(x_i), \ldots, \alpha^{k+s}(x_n)]=$

$-\epsilon(D_{\xi}+D_{\eta}, X_{i-1})[\alpha^{k+s}(x_1),  \ldots, [D_{\xi},D_{\eta}](x_i), \ldots, \alpha^{k+s}(x_n)].$

On the other hand, 
since ${\rm QC}(T)$ is a color $Hom$-Lie algebra, then, for arbitrary $x_1, \ldots, x_n \in T$, we have

$[[D_{\xi},D_{\eta}](x_1), \alpha^{k+s}(x_2), \ldots, \alpha^{k+s}(x_n)]=\epsilon(D_{\xi}+D_{\eta}, X_{i-1}) [\alpha^{k+s}(x_1),\ldots, [D_{\xi},D_{\eta}](x_i), \ldots, \alpha^{k+s}(x_n)]$.

Therefore
$$[\alpha^{k+s}(x_1),\ldots, [D_{\xi},D_{\eta}](x_i), \ldots, \alpha^{k+s}(x_n)]=0
\mbox{, hence }
[x_1',\ldots, [D_{\xi},D_{\eta}](x_i), \ldots, x_n']=0.$$

Hence, $[D_{\xi},D_{\eta}]=0.$




$(\Leftarrow)$\quad Is clear.\hfill$\Box$

\section{Quasiderivations of  multiplicative  $n$-ary $Hom$-$\Omega$ color algebras}

In this section we investigate the quasiderivations of the multiplicative $n$-ary $Hom$-$\Omega$ color algebra $T.$
We prove that $QDer(T)$ can be embedded in the derivation algebra of a larger multiplicative $n$-ary $Hom$-$\Omega$ color algebra $\breve{T}$ for the same variety of polynomial $Hom$-identities $\Omega.$
Moreover, we conclude  that $Der(\breve{T})$ has a direct sum decomposition when $Ann(T)=\{0\}.$

\


\begin{Lem}\label{33} 
Let $T$ be a multiplicative $n$-ary $Hom$-$\Omega$ color algebra over
       ${\mathbb F}$ and $t$ be an indeterminate. 
Define 
$$\breve{T}:= \{\Sigma(x\otimes t+y\otimes t^{n}) \ | \ x,y\in T\} \mbox{
and }\breve{\alpha}(\breve{T})= \{\Sigma ( \alpha(x) \otimes t + \alpha (y) \otimes t^n) \ | \  x,y\in T\}.$$
Endow $(\breve{T}, \breve{\alpha})$ with the multiplication
$$[x_1\otimes t^{i_1},x_2 \otimes
t^{i_2}, \ldots ,x_n\otimes t^{i_n}]=[x_1,x_2, \ldots, x_n]\otimes t^{\sum i_j},$$
for $ i_1, \ldots, i_n \in \{1,n\}$ (we put $t^k = 0$ if $k > n$).\\
Then $(\breve{T}, \breve{\alpha})$
is a multiplicative  $n$-ary $Hom$-$\Omega$ color algebra.
\end{Lem}

{\it Proof.}\quad Let the class of  $n$-ary $Hom$-$\Omega$ color algebras be defined by the  family $\{ f_k \}$ of color multilinear $Hom$-identities.
Then, for arbitrary $x_1, x_2, \ldots, x_m \in T$ and $i_j \in\{1,n\}$ we have
$$f_j(x_1\otimes t^{i_1},x_2\otimes t^{i_2}, \ldots, x_m\otimes t^{i_m})=
f_j(x_{1},x_{2}, \ldots, x_{m}) \otimes t^{\sum i_l}=0.$$

Therefore $\breve{T}$ is a $n$-ary $Hom$-$\Omega$ color algebra. \hfill$\Box$

For the sake of convenience, we write $xt ~ (xt^{n})$ instead of $x\otimes t ~ (x\otimes t^{n}).$

If $U$ is a $\mathbb{G}$-graded subspace of $T$ such that $T=U\oplus [T,\ldots,T],$ then
$$\breve{T}=Tt+Tt^{n}=Tt+Ut^{n}+[T,\ldots,T]t^{n}.$$
Now define a map $\varphi:{\rm QDer}(T)\rightarrow {\rm End}(\breve{T})$ by
$$\varphi(D)(at+ut^{n}+bt^{n})=D(a)t+D'(b)t^{n},$$
where $D\in {\rm QDer}($T$),$  $D'$ is a map related to $D$ by the definition of quasiderivation, $a\in
T,u\in U,b\in [T,\ldots,T]$.

\begin{Th}\label{34}
Let $T,\breve{T},\varphi$ be as above. Then

$(1)$ $\varphi$ is even; 

$(2)$ $\varphi$ is injective and $\varphi(D)$ does not depend on the choice of $D'$;

$(3)$ $\varphi({\rm QDer}(T))\subseteq {\rm Der}(\breve{T}).$

\end{Th}
{\it Proof.}\quad (1) Follows from the definition of $\varphi.$

(2) \quad If $\varphi(D_{1})=\varphi(D_{2}),$ then
for all $a\in T,b\in [T,\ldots ,T]$ and $u\in U$ we have
$$\varphi(D_{1})(at+ut^{n}+bt^{n})=\varphi(D_{2})(at+ut^{n}+bt^{n}),$$ or, in terms of $D_1, D_2,$
$$D_{1}(a)t+D'_{1}(b)t^{n}=D_{2}(a)t+D'_{2}(b)t^{n},$$ so $D_{1}(a)=D_{2}(a).$ 
Hence $D_{1}=D_{2},$ and $\varphi$ is injective.

Suppose that there exists $D''$ such that
$$\varphi(D)(at+ut^{n}+bt^{n})=D(a)t+D''(b)t^{n},$$ and
$$\sum  \epsilon(D, X_{i-1}) [\alpha^k(x_1), \ldots, D(x_i),\ldots, \alpha^k(x_n)]=D''([x_1, \ldots, x_n]),$$
then we have $$D'([x_1, \ldots, x_n])=D''([x_1, \ldots, x_n]),$$ thus $D'(b)=D''(b).$
Hence
$$\varphi(D)(at+ut^{n}+bt^{n})=D(a)t+D'(b)t^{n}=D(a)t+D''(b)t^{n},$$
which implies that $\varphi(D)$ is determined only by $D$.

(3)\quad We have $[x_1t^{i_1}, \ldots, x_n t^{i_n}]=[x_1 , \ldots, x_n ]t^{\sum i_j}=0$ for
all $\sum i_j \geq n+1$. Thus, to show that  $\varphi(D)\in {\rm
Der}(\breve{T}),$ we only need to check  that the following equality holds:
$$\varphi(D)([x_1t, \ldots, x_nt])=\sum \epsilon(D, X_{i-1}) [ \breve{\alpha}^k(x_1t),\ldots, \varphi(D)(x_it),\ldots, \breve{\alpha}^k(x_nt)].$$

For arbitrary $x_1, \ldots, x_n \in T$ we have
$$\varphi(D)([x_1t,\ldots, x_it, \ldots, x_nt])=\varphi(D)([x_1,\ldots, x_n]t^n)=
D'([x_1, \ldots, x_n])t^n$$
$$=\sum  \epsilon(D, X_{i-1}) [{\alpha}^k(x_1), \ldots, D(x_i),\ldots, {\alpha}^k(x_n)]t^n$$
$$=\sum  \epsilon(D, X_{i-1}) [{\alpha}^k(x_1t), \ldots, D(x_i)t,\ldots, {\alpha}^k(x_nt)]$$
$$=\sum  \epsilon(D, X_{i-1}) [\breve{\alpha}^k(x_1t), \ldots, \varphi(D)(x_it),\ldots, \breve{\alpha}^k(x_nt)].$$
Therefore, for all $D\in {\rm QDer}(T)$ we have $\varphi(D)\in {\rm
Der}(\breve{T})$. 
 \hfill$\Box$

\

\begin{Lem}\label{35}
Let $T$ be a multiplicative $n$-ary $Hom$-$\Omega$ color algebra such that
${\rm Ann}(T)=\{0\}$ and let $\breve{T},~\varphi$ be as previously defined.
Then $${\rm Der}(\breve{T})=\varphi({\rm QDer}(T))\oplus {\rm
ZDer}(\breve{T}).$$
\end{Lem}
{\it Proof.}\quad Since ${\rm Ann}(T)=\{0\}$, we have ${\rm
Ann}(\breve{T})=Tt^n.$
For $g\in {\rm Der}(\breve{T})$ we have
$g({\rm Ann}(\breve{T}))\subseteq {\rm Ann}(\breve{T}),$ hence
$g(Ut^n)\subseteq g({\rm Ann}(\breve{T}))\subseteq {\rm
Ann }(\breve{T})=Tt^n.$ Now define a map
$f:Tt+Ut^n+[T,\ldots,T]t^n\rightarrow Tt^n$ by
$$\ f(x)=\left\{\begin{array}{ll}g(x)\cap Tt^n,& x\in Tt ;\\
 g(x),& x\in Ut^n ;\\  0,& x\in [T,\ldots,T]t^n.\end{array}\right.$$

It is clear that $f$ is linear. Note that
$$f([\breve{T}, \ldots, \breve{T}])=f([T,\ldots ,T]t^n)=0,$$
$$~[\breve{\alpha}^k(\breve{T}), \ldots, f(\breve{T}), \ldots,  \breve{\alpha}^k(\breve{T})]
\subseteq [\alpha^k(T)t+{\alpha}^k(T)t^n,\ldots, Tt^n,\ldots, {\alpha}^k(T)t+{\alpha}^k(T)t^n]=0,$$
hence $f\in {\rm ZDer}(\breve{T}).$
Since $$(g-f)(Tt)=g(Tt)-g(Tt)\cap Tt^{n}=g(Tt)-Tt^{n}\subseteq Tt,~
(g-f)(Ut^n)=0,$$ and
$$(g-f)([T, \ldots, T]t^n)=g([\breve{T},\ldots ,\breve{T}])\subseteq
[\breve{T},\ldots, \breve{T}]=[T, \ldots, T]t^n,$$ there exist $D,~D'\in
{\rm End}(T)$ such that for all $a\in T,~b\in [T, \ldots, T]$,
$$(g-f)(at)=D(a)t,~ (g-f)(bt^n)=D'(b)t^n.$$ 
Since $g-f\in {\rm
Der}(\breve{T}),$ by the definition of ${\rm Der}(\breve{T})$ we
have
$$\sum \epsilon(g-f, A_{i-1}) [\breve{\alpha}^k(a_1t), \ldots, (g-f)(a_it), \ldots, \breve{\alpha}^k(a_nt)]=
(g-f)([a_1t, \ldots, a_nt]),$$
for all $a_1, \ldots, a_n\in T.$
Hence
$$\sum \epsilon(D, A_{i-1}) [{\alpha}^k(a_1), \ldots, D(a_i), \ldots, {\alpha}^k(a_n)]=D'([a_1,  \ldots, a_n]).$$ Thus
$D\in {\rm QDer}(T).$ 
Therefore, $g-f=\varphi(D)\in \varphi({\rm
QDer}(T))$, so ${\rm Der}(\breve{T})\subseteq \varphi({\rm
QDer}(T))+{\rm ZDer}(\breve{T}).$ By Lemma \ref{34}~$(3)$ we have
${\rm Der}(\breve{T})=\varphi({\rm QDer}(T))+{\rm ZDer}(\breve{T}).$

For any $f\in \varphi({\rm QDer}(T))\cap{\rm ZDer}(\breve{T})$
there exists an element $D\in {\rm QDer}(T)$ such that
$f=\varphi(D).$ Then
$$f(at+ut^n+bt^n)=\varphi(D)(at+ut^n+bt^n)=D(a)t+D'(b)t^n,$$ where $a\in T,b\in [T, \ldots, T].$

On the other hand, since $f\in {\rm ZDer}(\breve{T}),$ we have
$$f(at+bt^n+ut^n)\in {\rm Ann}(\breve{T})=Tt^n.$$ That is,
$D(a)=0,$ for all $a\in T$ and so $D=0.$ Hence $f=0.$

Therefore ${\rm Der}(\breve{T})=\varphi({\rm QDer}(T))\oplus {\rm
ZDer}(\breve{T})$ as desired. \hfill$\Box$

\end{document}